\documentclass[]{amsart}
\usepackage{amsmath, amsthm, amscd, amsfonts, amssymb, graphicx, color}
\usepackage[bookmarksnumbered, colorlinks, plainpages]{hyperref}

\theoremstyle{definition}

\theoremstyle{remark}

\numberwithin{equation}{section}

\begin{document}
\setcounter{page}{1}
\begin{center}
{\bf THE TOPOLOGICAL CENTERS  OF MODULE ACTIONS}
\end{center}

\title[]{}

\author[]{KAZEM HAGHNEJAD AZAR  }

\address{}

\dedicatory{}



\begin{abstract}
 In this article,  for Banach left and right module actions, we will extend some propositions from Lau and $\ddot{U}lger$ into general situations  and we establish  the relationships between  topological centers of module actions. We also introduce the new concepts as $Lw^*w$-property and $Rw^*w$-property for Banach $A-bimodule$ $B$ and we investigate the  relations between them and topological center of module actions. We have some applications in  dual groups.
 \end{abstract} \maketitle

\begin{center}
{\bf  1.Introduction and Preliminaries}
\end{center}

\noindent As is well-known [1], the second dual $A^{**}$ of $A$ endowed with the either Arens multiplications is a Banach algebra. The constructions of the two Arens multiplications in $A^{**}$ lead us to definition of topological centers for $A^{**}$ with respect both Arens multiplications. The topological centers of Banach algebras, module actions and applications of them  were introduced and discussed in [6, 8, 13, 14, 15, 16, 17, 21, 22], and they have attracted by some attentions.

\noindent Now we introduce some notations and definitions that we used
throughout  this paper.\\
 Let $A$ be a Banach algebra. We
say that a  net $(e_{\alpha})_{{\alpha}\in I}$ in $A$ is a left
approximate identity $(=LAI)$ [resp. right
approximate identity $(=RAI)$] if,
 for each $a\in A$,   $e_{\alpha}a\longrightarrow a$ [resp. $ae_{\alpha}\longrightarrow a$]. For $a\in A$
 and $a^\prime\in A^*$, we denote by $a^\prime a$
 and $a a^\prime$ respectively, the functionals on $A^*$ defined by $<a^\prime a,b>=<a^\prime,ab>=a^\prime(ab)$ and $<a a^\prime,b>=<a^\prime,ba>=a^\prime(ba)$ for all $b\in A$.
  The Banach algebra $A$ is embedded in its second dual via the identification
 $<a,a^\prime>$ - $<a^\prime,a>$ for every $a\in
A$ and $a^\prime\in
A^*$. We denote the set   $\{a^\prime a:~a\in A~ and ~a^\prime\in
  A^*\}$ and
  $\{a a^\prime:~a\in A ~and ~a^\prime\in A^*\}$ by $A^*A$ and $AA^*$, respectively, clearly these two sets are subsets of $A^*$.  Let $A$ has a $BAI$. If the
equality $A^*A=A^*,~~(AA^*=A^*)$ holds, then we say that $A^*$
factors on the left (right). If both equalities $A^*A=AA^*=A^*$
hold, then we say
that $A^*$  factors on both sides.
 Let $X,Y,Z$ be normed spaces and $m:X\times Y\rightarrow Z$ be a bounded bilinear mapping. Arens in [1] offers two natural extensions $m^{***}$ and $m^{t***t}$ of $m$ from $X^{**}\times Y^{**}$ into $Z^{**}$ as following:\\
 \noindent1. $m^*:Z^*\times X\rightarrow Y^*$,~~~~~given by~~~$<m^*(z^\prime,x),y>=<z^\prime, m(x,y)>$ ~where $x\in X$, $y\in Y$, $z^\prime\in Z^*$,\\
 2. $m^{**}:Y^{**}\times Z^{*}\rightarrow X^*$,~~given by $<m^{**}(y^{\prime\prime},z^\prime),x>=<y^{\prime\prime},m^*(z^\prime,x)>$ ~where $x\in X$, $y^{\prime\prime}\in Y^{**}$, $z^\prime\in Z^*$,\\
3. $m^{***}:X^{**}\times Y^{**}\rightarrow Z^{**}$,~ given by~ ~ ~$<m^{***}(x^{\prime\prime},y^{\prime\prime}),z^\prime>$ $=<x^{\prime\prime},m^{**}(y^{\prime\prime},z^\prime)>$ \\~where ~$x^{\prime\prime}\in X^{**}$, $y^{\prime\prime}\in Y^{**}$, $z^\prime\in Z^*$.\\
The mapping $m^{***}$ is the unique extension of $m$ such that $x^{\prime\prime}\rightarrow m^{***}(x^{\prime\prime},y^{\prime\prime})$ from $X^{**}$ into $Z^{**}$ is $weak^*-to-weak^*$ continuous for every $y^{\prime\prime}\in Y^{**}$, but the mapping $y^{\prime\prime}\rightarrow m^{***}(x^{\prime\prime},y^{\prime\prime})$ is not in general $weak^*-to-weak^*$ continuous from $Y^{**}$ into $Z^{**}$ unless $x^{\prime\prime}\in X$. Hence the first topological center of $m$ may  be defined as following
$$Z_1(m)=\{x^{\prime\prime}\in X^{**}:~~y^{\prime\prime}\rightarrow m^{***}(x^{\prime\prime},y^{\prime\prime})~~is~~weak^*-to-weak^*-continuous\}.$$
Let now $m^t:Y\times X\rightarrow Z$ be the transpose of $m$ defined by $m^t(y,x)=m(x,y)$ for every $x\in X$ and $y\in Y$. Then $m^t$ is a continuous bilinear map from $Y\times X$ to $Z$, and so it may be extended as above to $m^{t***}:Y^{**}\times X^{**}\rightarrow Z^{**}$.
 The mapping $m^{t***t}:X^{**}\times Y^{**}\rightarrow Z^{**}$ in general is not equal to $m^{***}$, see [1], if $m^{***}=m^{t***t}$, then $m$ is called Arens regular. The mapping $y^{\prime\prime}\rightarrow m^{t***t}(x^{\prime\prime},y^{\prime\prime})$ is $weak^*-to-weak^*$ continuous for every $y^{\prime\prime}\in Y^{**}$, but the mapping $x^{\prime\prime}\rightarrow m^{t***t}(x^{\prime\prime},y^{\prime\prime})$ from $X^{**}$ into $Z^{**}$ is not in general  $weak^*-to-weak^*$ continuous for every $y^{\prime\prime}\in Y^{**}$. So we define the second topological center of $m$ as
$$Z_2(m)=\{y^{\prime\prime}\in Y^{**}:~~x^{\prime\prime}\rightarrow m^{t***t}(x^{\prime\prime},y^{\prime\prime})~~is~~weak^*-to-weak^*-continuous\}.$$
It is clear that $m$ is Arens regular if and only if $Z_1(m)=X^{**}$ or $Z_2(m)=Y^{**}$. Arens regularity of $m$ is equivalent to the following
$$\lim_i\lim_j<z^\prime,m(x_i,y_j)>=\lim_j\lim_i<z^\prime,m(x_i,y_j)>,$$
whenever both limits exist for all bounded sequences $(x_i)_i\subseteq X$ , $(y_i)_i\subseteq Y$ and $z^\prime\in Z^*$, see [6, 18].\\
 The regularity of a normed algebra $A$ is defined to be the regularity of its algebra multiplication when considered as a bilinear mapping. Let $a^{\prime\prime}$ and $b^{\prime\prime}$ be elements of $A^{**}$, the second dual of $A$. By $Goldstin^,s$ Theorem [6, P.424-425], there are nets $(a_{\alpha})_{\alpha}$ and $(b_{\beta})_{\beta}$ in $A$ such that $a^{\prime\prime}=weak^*-\lim_{\alpha}a_{\alpha}$ ~and~  $b^{\prime\prime}=weak^*-\lim_{\beta}b_{\beta}$. So it is easy to see that for all $a^\prime\in A^*$,
$$\lim_{\alpha}\lim_{\beta}<a^\prime,m(a_{\alpha},b_{\beta})>=<a^{\prime\prime}b^{\prime\prime},a^\prime>$$ and
$$\lim_{\beta}\lim_{\alpha}<a^\prime,m(a_{\alpha},b_{\beta})>=<a^{\prime\prime}ob^{\prime\prime},a^\prime>,$$
where $a^{\prime\prime}b^{\prime\prime}$ and $a^{\prime\prime}ob^{\prime\prime}$ are the first and second Arens products of $A^{**}$, respectively, see [6, 14, 18].\\
The mapping $m$ is left strongly Arens irregular if $Z_1(m)=X$ and $m$ is right strongly Arens irregular if $Z_2(m)=Y$.\\
This paper is organized as follows.\\
a) In section two, for  a  Banach $A-bimodule$, we have
 \begin{enumerate}
 \item $a^{\prime\prime}\in Z_{B^{**}}(A^{**})$ if and only if $\pi_\ell^{****}(b^\prime,a^{\prime\prime})\in B^*$ for all $b^\prime\in B^*$.

\item $F\in Z_{B^{**}}((A^*A)^*)$ if and only if $\pi_\ell^{****}(g,F)\in B^{*}$ for all $g\in B^{*}$.
\item  $G\in Z_{(A^*A)^*}(B^{**})$ if and only if $\pi_r^{****}(g,G)\in A^*A$ for all $g\in B^{*}$.
\item  Let $B$ has a $BAI$ $(e_\alpha)_\alpha\subseteq A$ such that
$e_\alpha \stackrel{w^*} {\rightarrow}e^{\prime\prime}$. Then if ${Z}^t_{e^{**}}(B^{**})=B^{**}$ [ resp. ${Z}_{e^{**}}(B^{**})=B^{**}$] and $B^*$ factors on the left [resp. right], but not on the right [resp. left], then ${Z}_{B^{**}}(A^{**})\neq {Z}^t_{B^{**}}(A^{**})$.

\item $B^*A\subseteq wap_\ell(B)$
if and only if  $AA^{**}\subseteq Z_{B^{**}}(A^{**})$.
\item  Let $b^\prime\in B^*$. Then $b^\prime\in wap_\ell(B)$ if and only if the adjoint of  the mapping $\pi_\ell^*(b^\prime,~):A\rightarrow B^*$ is $weak^*-to-weak$ continuous.
\end{enumerate}
b) In section three, for  a Banach  $A-bimodule$ $B$, we define $Left-weak^*-to-weak$ property [=$Rw^*w-$ property]  and $Right-weak^*-to-weak$ property [=$Rw^*w-$ property] for Banach algebra $A$ and we show that
\begin{enumerate}
\item If $A^{**}=a_0A^{**}$ [resp.  $A^{**}=A^{**}a_0$] for some $a_0\in A$ and $a_0$ has $Rw^*w-$ property [resp. $Lw^*w-$ property], then $Z_{B^{**}}(A^{**})=A^{**}$.
\item If $B^{**}=a_0B^{**}$ [resp.  $B^{**}=B^{**}a_0$] for some $a_0\in A$ and $a_0$ has $Rw^*w-$ property [resp. $Lw^*w-$ property] with respect to $B$, then $Z_{A^{**}}(B^{**})=B^{**}$.
\item  If $B^*$ factors on the left [resp. right] with respect to $A$ and $A$ has $Rw^*w-$ property [resp. $Lw^*w-$ property], then $Z_{B^{**}}(A^{**})=A^{**}$.
\item  If $B^*$ factors on the left [resp. right] with respect to $A$ and $A$ has $Rw^*w-$ property [resp. $Lw^*w-$ property] with respect $B$, then $Z_{A^{**}}(B^{**})=B^{**}$.
\item  If $a_0\in A$ has $Rw^*w-$ property with respect to $B$, then $a_0A^{**}\subseteq Z_{B^{**}}(A^{**})$ and $a_0B^*\subseteq wap_\ell(B)$.
 \item    Assume that $AB^*\subseteq wap_\ell B$. If
$B^*$ strong factors on the left [resp. right], then $A$ has $Lw^*w-$ property [resp. $Rw^*w-$ property ] with respect to $B$.
\item Assume that $AB^*\subseteq wap_\ell B$. If
$B^*$ strong factors on the left [resp. right], then $A$ has $Lw^*w-$ property [resp. $Rw^*w-$ property ] with respect to $B$.\\\\
\end{enumerate}

\begin{center}
\textbf{{ 2. The topological centers of module actions}}
\end{center}

\noindent Let $B$ be a Banach $A-bimodule$, and let\\
$$\pi_\ell:~A\times B\rightarrow B~~~and~~~\pi_r:~B\times A\rightarrow B.$$
be the left and right module actions of $A$ on $B$. Then $B^{**}$ is a Banach $A^{**}-bimodule$ with module actions
$$\pi_\ell^{***}:~A^{**}\times B^{**}\rightarrow B^{**}~~~and~~~\pi_r^{***}:~B^{**}\times A^{**}\rightarrow B^{**}.$$
Similarly, $B^{**}$ is a Banach $A^{**}-bimodule$ with module actions\\
$$\pi_\ell^{t***t}:~A^{**}\times B^{**}\rightarrow B^{**}~~~and~~~\pi_r^{t***t}:~B^{**}\times A^{**}\rightarrow B^{**}.$$
We may therefore define the topological centers of the right and left module actions of $A$ on $B$ as follows:
$$Z_{A^{**}}(B^{**})=Z(\pi_r)=\{b^{\prime\prime}\in B^{**}:~the~map~~a^{\prime\prime}\rightarrow \pi_r^{***}(b^{\prime\prime}, a^{\prime\prime})~:~A^{**}\rightarrow B^{**}$$$$~is~~~weak^*-to-weak^*~continuous\}$$
$$Z_{B^{**}}(A^{**})=Z(\pi_\ell)=\{a^{\prime\prime}\in A^{**}:~the~map~~b^{\prime\prime}\rightarrow \pi_\ell^{***}(a^{\prime\prime}, b^{\prime\prime})~:~B^{**}\rightarrow B^{**}$$$$~is~~~weak^*-to-weak^*~continuous\}$$
$$Z_{A^{**}}^t(B^{**})=Z(\pi_\ell^t)=\{b^{\prime\prime}\in B^{**}:~the~map~~a^{\prime\prime}\rightarrow \pi_\ell^{t***}(b^{\prime\prime}, a^{\prime\prime})~:~A^{**}\rightarrow B^{**}$$$$~is~~~weak^*-to-weak^*~continuous\}$$
$$Z_{B^{**}}^t(A^{**})=Z(\pi_r^t)=\{a^{\prime\prime}\in A^{**}:~the~map~~b^{\prime\prime}\rightarrow \pi_r^{t***}(a^{\prime\prime}, b^{\prime\prime})~:~B^{**}\rightarrow B^{**}$$$$~is~~~weak^*-to-weak^*~continuous\}$$
We note also that if $B$ is a left(resp. right) Banach $A-module$ and $\pi_\ell:~A\times B\rightarrow B$~(resp. $\pi_r:~B\times A\rightarrow B$) is left (resp. right) module action of $A$ on $B$, then $B^*$ is a right (resp. left) Banach $A-module$. \\
We write $ab=\pi_\ell(a,b)$, $ba=\pi_r(b,a)$, $\pi_\ell(a_1a_2,b)=\pi_\ell(a_1,a_2b)$,  $\pi_r(b,a_1a_2)=\pi_r(ba_1,a_2)$,~\\
$\pi_\ell^*(a_1b^\prime, a_2)=\pi_\ell^*(b^\prime, a_2a_1)$,~
$\pi_r^*(b^\prime a, b)=\pi_r^*(b^\prime, ab)$,~ for all $a_1,a_2, a\in A$, $b\in B$ and  $b^\prime\in B^*$
when there is no confusion.\\

\noindent{\it{\bf Theorem 2-1.}} We have the following assertions.
\begin{enumerate}
\item  Assume that   $B$ is a Left Banach $A-module$. Then, $a^{\prime\prime}\in Z_{B^{**}}(A^{**})$ if and only if $\pi_\ell^{****}(b^\prime,a^{\prime\prime})\in B^*$ for all $b^\prime\in B^*$.
\item  Assume that   $B$ is a right Banach $A-module$. Then, $b^{\prime\prime}\in Z_{A^{**}}(B^{**})$ if and only if $\pi_r^{****}(b^\prime,b^{\prime\prime})\in A^*$ for all $b^\prime\in B^*$.
\end{enumerate}
\begin{proof}
\begin{enumerate}
\item Let $b^{\prime\prime}\in B^{**}$. Then, for every $a^{\prime\prime}\in Z_{B^{**}}(A^{**})$, we have
$$<\pi_\ell^{****}(b^\prime,a^{\prime\prime}),b^{\prime\prime}>=
<b^\prime,\pi_\ell^{***}(a^{\prime\prime},b^{\prime\prime})>=
<\pi_\ell^{***}(a^{\prime\prime},b^{\prime\prime}),b^\prime>$$
$$=<\pi_\ell^{t***t}(a^{\prime\prime},b^{\prime\prime}),b^\prime>=
<\pi_\ell^{t***}(b^{\prime\prime},a^{\prime\prime}),b^\prime>=
<b^{\prime\prime},\pi_\ell^{t**}(a^{\prime\prime},b^\prime)>.$$
It follow that $\pi_\ell^{****}(b^\prime,a^{\prime\prime})=\pi_\ell^{t**}(a^{\prime\prime},b^\prime)\in B^*$.\\
Conversely,  let $a^{\prime\prime}\in A^{**}$ and let $\pi_\ell^{****}(a^{\prime\prime},b^\prime)\in B^*$ for all $b^{\prime}\in B^{*}$.  Then for all $b^{\prime\prime}\in B^{**}$, we have
$$<\pi_\ell^{***}(a^{\prime\prime},b^{\prime\prime}),b^\prime>=
<b^\prime,\pi_\ell^{***}(a^{\prime\prime},b^{\prime\prime})>=<\pi_\ell^{****}(b^\prime,a^{\prime\prime}),
b^{\prime\prime}>$$
$$=<\pi_\ell^{t**}(a^{\prime\prime},b^\prime),b^{\prime\prime}>
=<b^{\prime\prime},\pi_\ell^{t**}(a^{\prime\prime},b^\prime)>=<\pi_\ell^{t***}(b^{\prime\prime},a^{\prime\prime}),
b^\prime>$$
$$=<\pi_\ell^{t***t}(a^{\prime\prime},b^{\prime\prime}),
b^\prime>.$$
Consequently $a^{\prime\prime}\in Z_{B^{**}}(A^{**})$.
\item Prof is similar to (1).
\end{enumerate}
\end{proof}
\noindent{\it{\bf Theorem 2-2.}} Assume that   $B$ is a Banach  $A-bimodule$. Then we have the following assertions.
\begin{enumerate}

\item $F\in Z_{B^{**}}((A^*A)^*)$ if and only if $\pi_\ell^{****}(g,F)\in B^{*}$ for all $g\in B^{*}$.
\item  $G\in Z_{(A^*A)^*}(B^{**})$ if and only if $\pi_r^{****}(g,G)\in A^*A$ for all $g\in B^{*}$.
\end{enumerate}
\begin{proof}
 \begin{enumerate}
 \item Let $F\in Z_{B^{**}}((A^*A)^*)$ and $(b_{\alpha}^{\prime\prime})_{\alpha}\subseteq B^{**}$ such that  $b^{\prime\prime}_{\alpha} \stackrel{w^*} {\rightarrow}b^{\prime\prime}$. Then for all $g\in B^{*}$, we have
$$<\pi_\ell^{****}(g,F),b^{\prime\prime}_{\alpha}>=<g,\pi_\ell^{***}(F,b^{\prime\prime}_{\alpha})>=
<\pi_\ell^{***}(F,b^{\prime\prime}_{\alpha}),g>$$
$$\rightarrow <\pi_\ell^{***}(F,b^{\prime\prime}),g>=<\pi_\ell^{****}(g,F),b^{\prime\prime}>.$$
Thus, we conclude that $\pi_\ell^{****}(g,F)\in (B^{**},weak^*)^*=B^{*}$.\\
Conversely, let $\pi_\ell^{****}(g,F)\in B^{*}$ for $F\in (A^*A)^*$ and $g\in B^{*}$. Assume that $b^{\prime\prime}\in B^{**}$ and  $(b_{\alpha}^{\prime\prime})_{\alpha}\subseteq B^{**}$ such that  $b^{\prime\prime}_{\alpha} \stackrel{w^*} {\rightarrow}b^{\prime\prime}$. Then
$$<\pi_\ell^{***}(F,b^{\prime\prime}_{\alpha}),g>=<g,\pi_\ell^{***}(F,b^{\prime\prime}_{\alpha})>
=<\pi_\ell^{****}(g,F),b^{\prime\prime}_{\alpha}>$$
$$=<b^{\prime\prime}_{\alpha},\pi_\ell^{****}(g,F)>\rightarrow <b^{\prime\prime},\pi_\ell^{****}(g,F)>=<\pi_\ell^{****}(g,F),b^{\prime\prime}>$$$$
=<\pi_\ell^{***}(F,b^{\prime\prime}),g>.$$
It follow that $F\in Z_{B^{**}}((A^*A)^*)$.
\item Proof is similar to (1).
 \end{enumerate}
 \end{proof}
In the proceeding theorems, if we take $B=A$, we obtain some parts of Lemma 3.1 from [14].\\

\noindent An element $e^{\prime\prime}$ of $A^{**}$ is said to be a mixed unit if $e^{\prime\prime}$ is a
right unit for the first Arens multiplication and a left unit for
the second Arens multiplication. That is, $e^{\prime\prime}$ is a mixed unit if
and only if,
for each $a^{\prime\prime}\in A^{**}$, $a^{\prime\prime}e^{\prime\prime}=e^{\prime\prime}o a^{\prime\prime}=a^{\prime\prime}$. By
[4, p.146], an element $e^{\prime\prime}$ of $A^{**}$  is  mixed
      unit if and only if it is a $weak^*$ cluster point of some BAI $(e_\alpha)_{\alpha \in I}$  in
      $A$.\\
Let $B$ be a Banach  $A-bimodule$ and $a^{\prime\prime}\in A^{**}$. We define the locally topological center of the left and right module actions of $a^{\prime\prime}$ on $B$, respectively, as follows\\
$$Z_{a^{\prime\prime}}^t(B^{**})=Z_{a^{\prime\prime}}^t(\pi_\ell^t)=\{b^{\prime\prime}\in B^{**}:~~~\pi^{t***t}_\ell(a^{\prime\prime},b^{\prime\prime})=
\pi^{***}_\ell(a^{\prime\prime},b^{\prime\prime})\},$$
$$Z_{a^{\prime\prime}}(B^{**})=Z_{a^{\prime\prime}}(\pi_r^t)=\{b^{\prime\prime}\in B^{**}:~~~\pi^{t***t}_r(b^{\prime\prime},a^{\prime\prime})=
\pi^{***}_r(b^{\prime\prime},a^{\prime\prime})\}.$$\\
Thus we have ~~~~~~~~~$$\bigcap_{a^{\prime\prime}\in A^{**}}Z_{a^{\prime\prime}}^t(B^{**})=Z_{A}^t(B^{**})=
Z(\pi_r^t),  $$ ~~~~~~~ ~~ $$~~~~~~~~~~~~~~~~~~~~~~~~~~~~\bigcap_{a^{\prime\prime}\in A^{**}}Z_{a^{\prime\prime}}(B^{**})=Z_{A}(B^{**})=
Z(\pi_r).$$\\\\
\noindent{\it{\bf Definition 2-3.}} Let $B$ be a left Banach  $A-module$ and $e^{\prime\prime}\in A^{**}$ be a mixed unit for $A^{**}$. We say that $e^{\prime\prime}$ is a left mixed unit for $B^{**}$, if
$$\pi_\ell^{***}(e^{\prime\prime},b^{\prime\prime})=\pi_\ell^{t***t}(e^{\prime\prime},b^{\prime\prime})=b^{\prime\prime},$$
for all $b^{\prime\prime}\in B^{**}$.\\
The definition of right mixed unit for $B^{**}$ is similar. $B^{**}$ has a mixed unit if it has left and right mixed unit that are equal.\\
It is clear that if   $e^{\prime\prime}\in A^{**}$ is a left (resp. right) unit for $B^{**}$ and $Z_{e^{\prime\prime}}(B^{**})=B^{**}$, then $e^{\prime\prime}$ is left (resp. right) mixed unit for $B^{**}$.\\

\noindent{\it{\bf Theorem 2-4.}} Let $B$ be a   Banach $A-bimodule$ with a $BAI$ $(e_\alpha)_\alpha$ such that
$e_\alpha \stackrel{w^*} {\rightarrow}e^{\prime\prime}$. Then if ${Z}^t_{e^{**}}(B^{**})=B^{**}$ [ resp. ${Z}_{e^{**}}(B^{**})=B^{**}$] and $B^*$ factors on the left [resp. right], but not on the right [resp. left], then ${Z}_{B^{**}}(A^{**})\neq {Z}^t_{B^{**}}(A^{**})$.
\begin{proof} Suppose that  $B^*$ factors on the left with respect to $A$, but not on the right. Let $(e_{\alpha})_{\alpha}\subseteq A$ be a BAI for $A$ such that  $e_{\alpha} \stackrel{w^*} {\rightarrow}e^{\prime\prime}$. Thus for all $b^{\prime}\in B^{*}$ there are $a\in A$ and $x^\prime\in B^*$ such that $x^\prime a=b^\prime$. Then for all $b^{\prime\prime}\in B^{**}$ we have
$$<\pi_\ell^{***}(e^{\prime\prime},b^{\prime\prime}),b^\prime>
=<e^{\prime\prime},\pi_\ell^{**}(b^{\prime\prime},b^\prime)>=
\lim_\alpha<\pi_\ell^{**}(b^{\prime\prime},b^\prime),e_{\alpha}>$$
$$=\lim_\alpha<b^{\prime\prime},\pi_\ell^{*}(b^\prime,e_{\alpha})>
=\lim_\alpha<b^{\prime\prime},\pi_\ell^{*}(x^\prime a,e_{\alpha})>$$
$$=\lim_\alpha<b^{\prime\prime},\pi_\ell^{*}(x^\prime ,ae_{\alpha})>
=\lim_\alpha<\pi_\ell^{**}(b^{\prime\prime},x^\prime) ,ae_{\alpha}>$$$$=<\pi_\ell^{**}(b^{\prime\prime},x^\prime) ,a>
=<b^{\prime\prime},b^{\prime}>.$$
Thus $\pi^{***}_\ell(e^{\prime\prime},b^{\prime\prime})=b^{\prime\prime}$ consequently $B^{**}$ has left unit $A^{**}-module$. It follows that
$e^{\prime\prime}\in {Z}_{B^{**}}(A^{**})$.  If we take
${Z}_{B^{**}}(A^{**})= {Z}^t_{B^{**}}(A^{**})$, then $e^{\prime\prime}\in {Z}^t_{B^{**}}(A^{**})$. Then the mapping $b^{\prime\prime}\rightarrow\pi_r^{t***t}(b^{\prime\prime},e^{\prime\prime})$ is $weak^*-to-weak^*$ continuous from $B^{**}$ into $B^{**}$. Since $e_\alpha \stackrel{w^*} {\rightarrow}e^{\prime\prime}$,
$\pi_r^{t***t}(b^{\prime\prime},e_\alpha)\stackrel{w^*} {\rightarrow}\pi_r^{t***t}(b^{\prime\prime},e^{\prime\prime})$. Let $b^\prime\in B^*$ and $(b_\beta)_\beta\subseteq B$ such that $b_\beta \stackrel{w^*} {\rightarrow}b^{\prime\prime}$. Since ${Z}^t_{e^{**}}(B^{**})=B^{**}$, we have the following quality
$$<\pi_r^{t***t}(b^{\prime\prime},e^{\prime\prime}), b^\prime>=\lim_\alpha<\pi_r^{t***t}(b^{\prime\prime},e_\alpha), b^\prime>=\lim_\alpha<\pi_r^{t***}(e_\alpha,b^{\prime\prime}), b^\prime>$$$$=\lim_\alpha\lim_\beta<\pi_r^{t***}(e_\alpha,b_\beta), b^\prime>
=
\lim_\alpha\lim_\beta<\pi_r^{}(b_\beta ,e_\alpha), b^\prime>$$$$=\lim_\alpha\lim_\beta< b^\prime, \pi_r^{}(b_\beta ,e_\alpha)>
=\lim_\beta\lim_\alpha< b^\prime, \pi_r^{}(b_\beta ,e_\alpha)>$$$$=\lim_\beta< b^\prime, b_\beta >=
<b^{\prime\prime},b^\prime>.$$
 Thus $\pi_r^{t***t}(b^{\prime\prime},e^{\prime\prime})=\pi_r^{***}(b^{\prime\prime},e^{\prime\prime})=b^{\prime\prime}$.
 It follows that $B^{\prime\prime}$ has a right unit.
 Suppose that $b^{\prime\prime}\in B^{**}$ and $(b_\beta)_\beta\subseteq B$ such that $b_\beta \stackrel{w^*} {\rightarrow}b^{\prime\prime}$.  Then for all $b^\prime\in B^*$ we have
$$<b^{\prime\prime},b^{\prime}>=<\pi_r^{***}(b^{\prime\prime},e^{\prime\prime}),b^\prime>=
<b^{\prime\prime},\pi_r^{**}(e^{\prime\prime},b^\prime)>=\lim_\beta<\pi_r^{**}(e^{\prime\prime},b^\prime),b_\beta>$$
$$=\lim_\beta <e^{\prime\prime},\pi_r^{*}(b^\prime,b_\beta)>
=\lim_\beta \lim_\alpha<\pi_r^{*}(b^\prime,b_\beta),e_\alpha>$$
$$=\lim_\beta \lim_\alpha<\pi_r^{*}(b^\prime,b_\beta),e_\alpha>=\lim_\beta \lim_\alpha<b^\prime,\pi_r(b_\beta,e_\alpha)>$$
$$=\lim_\alpha \lim_\beta<\pi_r^{***}(b_\beta,e_\alpha),b^\prime>
=\lim_\alpha \lim_\beta<b_\beta,\pi_r^{**}(e_\alpha,b^\prime)>$$
$$= \lim_\alpha<b^{\prime\prime},\pi_r^{**}(e_\alpha,b^\prime)>.$$
It follows that $weak-\lim_\alpha\pi_r^{**}(e_\alpha,b^\prime)=b^\prime$. So by Cohen Factorization Theorem, $B^*$ factors on the right that is contradiction.

\end{proof}

\noindent{\it{\bf Corollary 2-5.}} Let $B$ be a   Banach $A-bimodule$ and $e^{\prime\prime}\in A^{**}$ be a left mixed unit for $B^{**}$. If $B^*$ factors on the left, but not on the right, then ${Z}_{B^{**}}(A^{**})\neq {Z}^t_{B^{**}}(A^{**})$.\\

\noindent In the proceeding corollary, if we take $B=A$, then it is clear ${Z}^t_{e^{**}}(A^{**})=A^{**}$, and so we obtain Proposition 2.10 from [14].\\

\noindent{\it{\bf Theorem 2-6.}} Suppose that $B$ is a weakly complete Banach space. Then we have the following assertions.
\begin{enumerate}
\item Let $B$ be a Left Banach $A-module$ and $e^{\prime\prime}$ be a left mixed unit for $B^{**}$. If $AB^{**}\subseteq B$, then $B$ is reflexive.
\item Let $B$ be a right Banach $A-module$ and $e^{\prime\prime}$ be a right mixed unit for $B^{**}$. If $Z_{A^{**}}(B^{**})A\subseteq B$, then $Z_{A^{**}}(B^{**})= B$.
\end{enumerate}
\begin{proof}
\begin{enumerate}
\item Assume that $b^{\prime\prime}\in B^{**}$. Since $e^{\prime\prime}$ is also mixed unit for $A^{**}$, there is a $BAI$  $(e_{\alpha})_{\alpha}\subseteq A$ for $A$ such that $e_{\alpha} \stackrel{w^*} {\rightarrow}e^{\prime\prime}$. Then
$\pi_\ell^{***}(e_{\alpha},b^{\prime\prime}) \stackrel{w^*} {\rightarrow}\pi_\ell^{***}(e^{\prime\prime},b^{\prime\prime})=b^{\prime\prime}$ in $B^{**}$. Since $AB^{**}\subseteq B$, we have
$\pi_\ell^{***}(e_{\alpha},b^{\prime\prime})\in B$.  Consequently $\pi_\ell^{***}(e_{\alpha},b^{\prime\prime}) \stackrel{w} {\rightarrow}\pi_\ell^{***}(e^{\prime\prime},b^{\prime\prime})=b^{\prime\prime}$ in $B$.
Since $B$ is a weakly complete, $b^{\prime\prime}\in B$, and so $B$ is reflexive.
\item  Since  $b^{\prime\prime}\in Z_{A^{**}}(B^{**})$, we have  $\pi_r^{***}(b^{\prime\prime},e_{\alpha}) \stackrel{w^*} {\rightarrow}\pi_r^{***}(b^{\prime\prime},e^{\prime\prime})=b^{\prime\prime}$ in $B^{**}$.
    Since $Z_{A^{**}}(B^{**})A\subseteq B$, $\pi_r^{***}(b^{\prime\prime},e_{\alpha})\in B$. Consequently we have  $\pi_r^{***}(b^{\prime\prime},e_{\alpha}) \stackrel{w} {\rightarrow}\pi_r^{***}(b^{\prime\prime},e^{\prime\prime})=b^{\prime\prime}$ in $B$. It follows that $b^{\prime\prime}\in B$, since $B$ is a weakly complete.
\end{enumerate}
\end{proof}

\noindent A functional $a^\prime$ in $A^*$ is said to be $wap$ (weakly almost
 periodic) on $A$ if the mapping $a\rightarrow a^\prime a$ from $A$ into
 $A^{*}$ is weakly compact. The procceding definition to the equivalent following condition, see [6, 14, 18].\\
 For any two net $(a_{\alpha})_{\alpha}$ and $(b_{\beta})_{\beta}$
 in $\{a\in A:~\parallel a\parallel\leq 1\}$, we have\\
$$\\lim_{\alpha}\\lim_{\beta}<a^\prime,a_{\alpha}b_{\beta}>=\\lim_{\beta}\\lim_{\alpha}<a^\prime,a_{\alpha}b_{\beta}>,$$
whenever both iterated limits exist. The collection of all $wap$
functionals on $A$ is denoted by $wap(A)$. Also we have
$a^{\prime}\in wap(A)$ if and only if $<a^{\prime\prime}b^{\prime\prime},a^\prime>=<a^{\prime\prime}ob^{\prime\prime},a^\prime>$ for every $a^{\prime\prime},~b^{\prime\prime} \in
A^{**}$.\\

\noindent{\it{\bf Definition 2-7.}} Let $B$ be a left Banach $A-module$. Then, $b^\prime\in B^*$ is said to be left weakly almost periodic functional if the set $\{\pi_\ell(b^\prime,a):~a\in A,~\parallel a\parallel\leq 1\}$ is relatively weakly compact. We denote by $wap_\ell(B)$ the closed subspace of $B^*$ consisting of all the left weakly almost periodic functionals in $B^*$.\\
The definition of the right weakly almost periodic functional ($=wap_r(B)$) is the same.\\
By  [18], the definition of $wap_\ell(B)$ is equivalent to the following $$<\pi_\ell^{***}(a^{\prime\prime},b^{\prime\prime}),b^\prime>=
<\pi_\ell^{t***t}(a^{\prime\prime},b^{\prime\prime}),b^\prime>$$
for all $a^{\prime\prime}\in A^{**}$ and $b^{\prime\prime}\in B^{**}$.
Thus, we can write \\
$$wap_\ell(B)=\{ b^\prime\in B^*:~<\pi_\ell^{***}(a^{\prime\prime},b^{\prime\prime}),b^\prime>=
<\pi_\ell^{t***t}(a^{\prime\prime},b^{\prime\prime}),b^\prime>~~$$$$for~~all~~a^{\prime\prime}\in A^{**},~b^{\prime\prime}\in B^{**}\}.$$\\
\noindent{\it{\bf Theorem 2-8.}} Suppose that $B$ is a left Banach $A-module$. Consider the following statements.
\begin{enumerate}

\item $B^*A\subseteq wap_\ell(B)$.
\item $AA^{**}\subseteq Z_{B^{**}}(A^{**})$.
\item $AA^{**}\subseteq AZ_{B^{**}}((A^{*}A)^*)$.
\end{enumerate}
Then, we have $(1)\Leftrightarrow (2)\Leftarrow (3)$.

\begin{proof} $(1)\Rightarrow (2)$\\
 Let $(b_{\alpha}^{\prime\prime})_{\alpha}\subseteq B^{**}$ such that  $b^{\prime\prime}_{\alpha} \stackrel{w^*} {\rightarrow}b^{\prime\prime}$. Then for all $a\in A$ and $a^{\prime\prime}\in A^{**}$, we have
$$<\pi_\ell^{***}(aa^{\prime\prime},b^{\prime\prime}_{\alpha}),b^\prime>=
<aa^{\prime\prime},\pi_\ell^{**}(b^{\prime\prime}_{\alpha},b^\prime)>=
<a^{\prime\prime},\pi_\ell^{**}(b^{\prime\prime}_{\alpha},b^\prime)a>$$
$$=<a^{\prime\prime},\pi_\ell^{**}(b^{\prime\prime}_{\alpha},b^\prime a)>
=<\pi_\ell^{***}(a^{\prime\prime},b^{\prime\prime}_{\alpha}),b^\prime a)>
\rightarrow<\pi_\ell^{***}(a^{\prime\prime},b^{\prime\prime}),b^\prime a)>$$
$$=<\pi_\ell^{***}(aa^{\prime\prime},b^{\prime\prime}),b^\prime )>.$$
Hence $aa^{\prime\prime}\in Z_{B^{**}}(A^{**})$.\\
$(2)\Rightarrow (1)$\\
Let $a\in A$ and $b^\prime\in B^*$. Then
$$<\pi_\ell^{***}(a^{\prime\prime},b^{\prime\prime}_{\alpha}),b^\prime a>=
<a\pi_\ell^{***}(a^{\prime\prime},b^{\prime\prime}_{\alpha}),b^\prime >=
<\pi_\ell^{***}(aa^{\prime\prime},b^{\prime\prime}_{\alpha}),b^\prime >$$
$$=<\pi_\ell^{t***t}(aa^{\prime\prime},b^{\prime\prime}_{\alpha}),b^\prime >
=<\pi_\ell^{t***t}(a^{\prime\prime},b^{\prime\prime}_{\alpha}),b^\prime a>.$$
It follow that $b^\prime a\in wap_\ell(B)$.\\
$(3)\Rightarrow (2)$\\
Since $AZ_{B^{**}}((A^{*}A)^*)\subseteq Z_{B^{**}}(A^{**})$,  proof is hold.
\end{proof}
\noindent In the proceeding theorem, if we take $B=A$, then we obtain Theorem 3.6 from [14] and  the same as proceeding theorem, we can claim the following assertions:\\
If $B$ is a right Banach $A-module$, then for the following statements we have\\
 $(1)\Leftrightarrow (2)\Leftarrow (3)$.
 \begin{enumerate}

 \item $AB^*\subseteq wap_r(B)$.
\item $A^{**}A\subseteq Z_{B^{**}}(A^{**})$.
\item $A^{**}A\subseteq Z_{B^{**}}((A^{*}A)^*)A$.
\end{enumerate}
The proof of the this assertion is  similar to  proof of Theorem 2-8.\\\\
{\it{\bf Corollary 2-9.}} Suppose that $B$ is a  Banach $A-bimodule$. Then if $A$ is a left [resp. right] ideal in $A^{**}$, then $B^*A\subseteq wap_\ell(B)$ [resp. $AB^*\subseteq wap_r(B)$].\\\\
\noindent{\it{\bf Example 2-10.}} Suppose that $1\leq p\leq\infty$ and $q$ is conjugate of $p$. We know that if $G$ is compact, then $L^1(G)$ is a two-sided ideal in its second dual of it. By proceeding Theorem we have  $L^q(G)*L^1(G)\subseteq wap_\ell(L^p(G))$ and $L^1(G)*L^q(G)\subseteq wap_r(L^p(G))$.\\
Also if $G$ is finite, then $L^q(G)\subseteq wap_\ell(L^p(G))\cap wap_r(L^p(G))$. Hence we conclude that\\
$$Z_{L^1(G)^{**}}(L^p(G)^{**})=L^p(G)~~~~ and~~Z_{L^p(G)^{**}}(L^1(G)^{**})=L^1(G).$$\\
\noindent{\it{\bf Theorem 2-11.}} We have the following assertions.
\begin{enumerate}
\item Suppose that $B$ is a left Banach $A-module$ and $b^\prime\in B^*$. Then $b^\prime\in wap_\ell(B)$ if and only if the adjoint of  the mapping $\pi_\ell^*(b^\prime,~):A\rightarrow B^*$ is $weak^*-to-weak$ continuous.
\item Suppose that $B$ is a right Banach $A-module$ and $b^\prime\in B^*$. Then $b^\prime\in wap_r(B)$ if and only if the adjoint of the mapping $\pi_r^*(b^\prime,~):B\rightarrow A^*$ is $weak^*-to-weak$ continuous.
\end{enumerate}
\begin{proof}
\begin{enumerate}
 \item Assume that $b^\prime\in wap_\ell(B)$ and $\pi_\ell^*(b^\prime,~)^*:B^{**}\rightarrow A^*$ is the adjoint of $\pi_\ell^*(b^\prime,~)$. Then for every $b^{\prime\prime}\in B^{**}$ and $a\in A$, we have
$$<\pi_\ell^*(b^\prime,~)^*b^{\prime\prime},a>=<b^{\prime\prime},\pi_\ell^*(b^\prime,a)>.$$
Suppose $(b_{\alpha}^{\prime\prime})_{\alpha}\subseteq B^{**}$ such that  $b^{\prime\prime}_{\alpha} \stackrel{w^*} {\rightarrow}b^{\prime\prime}$ and $a^{\prime\prime}\in A^{**}$ and $(a_{\beta})_{\beta}\subseteq A$ such that
 $a_{\beta}\stackrel{w^*} {\rightarrow}a^{\prime\prime}$.
 By easy calculation, for all $y^{\prime\prime}\in B^{**}$ and $y^\prime \in B^*$,  we have
  $$<\pi_\ell^*(y^\prime,~)^*,y^{\prime\prime}>= \pi_\ell^{**}(y^{\prime\prime},y^\prime).$$
 Since $b^\prime\in wap_\ell(B)$,
$$<\pi_\ell^{***}(a^{\prime\prime},b_{\alpha}^{\prime\prime}),b^{\prime}>\rightarrow <\pi_\ell^{***}(a^{\prime\prime},b^{\prime\prime}),b^{\prime}>.$$
Then we have the following statements
$$\lim_\alpha<a^{\prime\prime},\pi_\ell^*(b^\prime,~)^*b_{\alpha}^{\prime\prime}>=
\lim_\alpha<a^{\prime\prime},\pi_\ell^{**}(b_{\alpha}^{\prime\prime},b^\prime)>$$
$$=
\lim_\alpha<\pi_\ell^{***}(a^{\prime\prime},b_{\alpha}^{\prime\prime}),b^\prime>=
<\pi_\ell^{***}(a^{\prime\prime},b^{\prime\prime}),b^\prime>$$
$$=<a^{\prime\prime},\pi_\ell^*(b^\prime,~)^*b^{\prime\prime}>.$$
It follow that the adjoint of the mapping $\pi_\ell^*(b^\prime,~):A\rightarrow B^*$ is $weak^*-to-weak$ continuous.\\
Conversely, let the adjoint of the mapping $\pi_\ell^*(b^\prime,~):A\rightarrow B^*$ is $weak^*-to-weak$ continuous.
Suppose $(b_{\alpha}^{\prime\prime})_{\alpha}\subseteq B^{**}$ such that  $b^{\prime\prime}_{\alpha} \stackrel{w^*} {\rightarrow}b^{\prime\prime}$ and $b^\prime \in B^*$. Then for every $a^{\prime\prime}\in A^{**}$, we have
$$\lim_\alpha<\pi_\ell^{***}(a^{\prime\prime},b_{\alpha}^{\prime\prime}),b^\prime>=
\lim_\alpha<a^{\prime\prime},\pi_\ell^{**}(b_{\alpha}^{\prime\prime},b^\prime)>$$$$=
\lim_\alpha<a^{\prime\prime},\pi_\ell^*(b^\prime,~)^*b_{\alpha}^{\prime\prime}>=
<a^{\prime\prime},\pi_\ell^*(b^\prime,~)^*b^{\prime\prime}>=
<\pi_\ell^{***}(a^{\prime\prime},b^{\prime\prime}),b^\prime>.$$
It follow that $b^\prime\in wap_\ell(B)$.\\
\item proof is similar to (1).
\end{enumerate}
\end{proof}

\noindent {\it{\bf Corollary 2-12.}} Let $A$ be a Banach algebra.
 Assume that $a^\prime\in A^*$ and $T_{a^\prime}$ is the linear operator from
$A$ into $A^*$ defined by $T_{a^\prime} a=a^\prime a$. Then, $a^\prime\in wap(A)$ if and
only if the adjoint of $T_{a^\prime}$ is
$weak^*-to-weak$ continuous. So $A$ is Arens regular if and only if the adjoint of the mapping  $T_{a^\prime} a=a^\prime a$ is $weak^*-to-weak$ continuous for every $a^\prime\in A^*$. \\

\begin{center}
\textbf{{ 3. $Lw^*w$-property and $Rw^*w$-property}}
\end{center}

\noindent In this section,  we introduce the new definition as $Left-weak^*-to-weak$ property and $Right-weak^*-to-weak$ property for Banach algebra $A$ and make some relations between these concepts and topological centers of module actions. As some conclusion, we have $Z_{L^1(G)^{**}}{(M(G)^{**})}\neq {M(G)^{**}}$ where $G$ is a locally compact group. If $G$ is finite, we have $Z_{M(G)^{**}}{(L^1(G)^{**})}= {L^1(G)^{**}}$ and
$Z_{L^1(G)^{**}}{(M(G)^{**})}= {M(G)^{**}}$.\\

\noindent{\it{\bf Definition 3-1.}} Let $B$ be a left Banach $A-module$. We say that $a\in A$ has $Left-weak^*-to-weak$ property ($=Lw^*w-$ property) with respect to $B$, if for all $(b_{\alpha})_{\alpha}\subseteq B^*$, $ab_\alpha^\prime\stackrel{w^*} {\rightarrow}0$ implies $ab_\alpha^\prime\stackrel{w} {\rightarrow}0$. If every $a\in A$ has $Lw^*w-$ property with respect to $B$, then we say that $A$ has $Lw^*w-$ property with respect to $B$. The definition of the $Right-weak^*-to-weak$ property ($=Rw^*w-$ property) is the same.\\
We say that $a\in A$ has $weak^*-to-weak$  property ($=w^*w-$ property) with respect to $B$ if it has $Lw^*w-$ property and $Rw^*w-$ property with respect to $B$.\\
If $a\in A$ has $Lw^*w-$ property with respect to itself, then we say that $a\in A$ has $Lw^*w-$ property.\\\\
For proceeding definition, we have  some examples and remarks as follows.\\
a) If $B$ is Banach $A$-bimodule and reflexive, then $A$ has $w^*w-$property with respect to $B$. Then\\
i) $L^1(G)$, $M(G)$ and $A(G)$ have $w^*w-$property when $G$ is finite.\\
ii) Let $G$ be locally compact group.  $L^1(G)$ [resp. $M(G)$] has $w^*w-$property [resp. $Lw^*w-$ property ]with respect to $L^p(G)$ whenever $p>1$.\\
b) Suppose that $B$ is a left Banach $A-module$ and $e$ is left unit element of $A$ such that $eb=b$ for all $b\in B$. If $e$ has $Lw^*w-$ property, then $B$ is reflexive.\\
c) If $S$ is a compact semigroup, then $C^+(S)=\{f\in C(S):~f>0\}$ has $w^*w-$property.\\\\\\
\noindent{\it{\bf Theorem 3-2.}} Suppose that $B$ is a Banach  $A-bimodule$. Then we have the following assertions.
 \begin{enumerate}
\item If $A^{**}=a_0A^{**}$ [resp.  $A^{**}=A^{**}a_0$] for some $a_0\in A$ and $a_0$ has $Rw^*w-$ property [resp. $Lw^*w-$ property], then $Z_{B^{**}}(A^{**})=A^{**}$.
\item If $B^{**}=a_0B^{**}$ [resp.  $B^{**}=B^{**}a_0$] for some $a_0\in A$ and $a_0$ has $Rw^*w-$ property [resp. $Lw^*w-$ property] with respect to $B$, then $Z_{A^{**}}(B^{**})=B^{**}$.
\end{enumerate}
\begin{proof}
\begin{enumerate}
\item Suppose that $A^{**}=a_0A^{**}$  for some $a_0\in A$ and $a_0$ has $Rw^*w-$ property. Let   $(b_{\alpha}^{\prime\prime})_{\alpha}\subseteq B^{**}$ such that  $b^{\prime\prime}_{\alpha} \stackrel{w^*} {\rightarrow}b^{\prime\prime}$. Then for all $a\in A$ and $b^\prime\in B^{*}$, we have
$$<\pi_\ell^{**}(b_\alpha^{\prime\prime},b^{\prime}), a>=<b_\alpha^{\prime\prime},\pi_\ell^{*}(b^{\prime}, a)>\rightarrow
<b^{\prime\prime},\pi_\ell^{*}(b^{\prime}, a)>=<\pi_\ell^{**}(b^{\prime\prime},b^{\prime}), a>,$$
it follow that $\pi_\ell^{**}(b_\alpha^{\prime\prime},b^{\prime})\stackrel{w^*} {\rightarrow}\pi_\ell^{**}(b^{\prime\prime},b^{\prime})$. Also we can write  $\pi_\ell^{**}(b_\alpha^{\prime\prime},b^{\prime})a_0\stackrel{w^*} {\rightarrow}\pi_\ell^{**}(b^{\prime\prime},b^{\prime})a_0$. Since $a_0$ has $Rw^*w-$ property, $\pi_\ell^{**}(b_\alpha^{\prime\prime},b^{\prime})a_0\stackrel{w} {\rightarrow}\pi_\ell^{**}(b^{\prime\prime},b^{\prime})a_0$. Now let $a^{\prime\prime}\in A^{**}$. Then there is
$x^{\prime\prime}\in A^{**}$ such that  $a^{\prime\prime}=a_0x^{\prime\prime}$ consequently we have
$$<\pi_\ell^{***}(a^{\prime\prime},b^{\prime\prime}_{\alpha}),b^\prime>=
<a^{\prime\prime},\pi_\ell^{**}(b^{\prime\prime}_{\alpha},b^\prime)>=
<x^{\prime\prime},\pi_\ell^{**}(b^{\prime\prime}_{\alpha},b^\prime)a_0>$$
$$\rightarrow<x^{\prime\prime},\pi_\ell^{**}(b^{\prime\prime},b^\prime)a_0>=
<\pi_\ell^{***}(a^{\prime\prime},b^{\prime\prime}_{\alpha}),b^\prime>.$$
We conclude that $a^{\prime\prime}\in Z_{B^{**}}(A^{**})$. Proof of the next part is the same as the proceeding proof.\\
\item Let $B^{**}=a_0B^{**}$  for some $a_0\in A$ and $a_0$ has $Rw^*w-$ property  with respect to $B$. Assume that
$(a_{\alpha}^{\prime\prime})_{\alpha}\subseteq A^{**}$ such that  $a^{\prime\prime}_{\alpha} \stackrel{w^*} {\rightarrow}a^{\prime\prime}$. Then for all $b\in B$, we have
$$<\pi_r^{**}(a_\alpha^{\prime\prime},b^{\prime}), b>=<a_\alpha^{\prime\prime},\pi_r^{**}(b^{\prime}, b)>
\rightarrow<a^{\prime\prime},\pi_r^{**}(b^{\prime}, b)>=<\pi_r^{**}(a^{\prime\prime},b^{\prime}), b>.$$
We conclude that $\pi_r^{**}(a_\alpha^{\prime\prime},b^{\prime})\stackrel{w^*} {\rightarrow}
\pi_r^{**}(a^{\prime\prime},b^{\prime})$ then we have $\pi_r^{**}(a_\alpha^{\prime\prime},b^{\prime})a_0\stackrel{w^*} {\rightarrow}
\pi_r^{**}(a^{\prime\prime},b^{\prime})a_0$. Since $a_0$ has $Rw^*w-$ property with respect to $B$,
$\pi_r^{**}(a_\alpha^{\prime\prime},b^{\prime})a_0\stackrel{w} {\rightarrow}
\pi_r^{**}(a^{\prime\prime},b^{\prime})a_0$.\\
Now let $b^{\prime\prime}\in B^{**}$. Then there is $x^{\prime\prime}\in B^{**}$ such that $b^{\prime\prime}=
a_0x^{\prime\prime}$. Hence, we have
$$<\pi_r^{***}(b^{\prime\prime},a_\alpha^{\prime\prime}), b^\prime>=<b^{\prime\prime},\pi_r^{**}(a^{\prime\prime}_\alpha, b^\prime)>=<a_0x^{\prime\prime},\pi_r^{**}(a^{\prime\prime}_\alpha, b^\prime)>$$
$$=<x^{\prime\prime},\pi_r^{**}(a^{\prime\prime}_\alpha, b^\prime)a_0>\rightarrow <x^{\prime\prime},\pi_r^{**}(a^{\prime\prime}, b^\prime)a_0>=<b^{\prime\prime},\pi_r^{**}(a^{\prime\prime}, b^\prime)>$$
$$=<\pi_r^{***}(b^{\prime\prime},a^{\prime\prime}), b^\prime>.$$
It follow that $b^{\prime\prime}\in Z_{A^{**}}(B^{**})$. The next part is similar to the proceeding proof.
\end{enumerate}
\end{proof}
\noindent{\it{\bf Example 3-3.}} i) Let $G$ be a locally compact group. Since $M(G)$ is a Banach $L^1(G)$-bimodule and the unit element of  $M(G)$ has not $Lw^*w-$ property or $Rw^*w-$ property, by Theorem 2-3, $Z_{L^1(G)^{**}}{(M(G)^{**})}\neq {M(G)^{**}}$.\\
ii) If $G$ is finite, then by Theorem 2-3, we have $Z_{M(G)^{**}}{(L^1(G)^{**})}= {L^1(G)^{**}}$ and
$Z_{L^1(G)^{**}}{(M(G)^{**})}= {M(G)^{**}}$.\\\\
Assume that  $B$ is a Banach $A-bimodule$. We say that  $B$ factors on the left (right) with respect to $A$ if $B=BA~(B=AB)$. We say that $B$ factors on both sides, if $B=BA=AB$.\\

\noindent{\it{\bf Theorem 3-4.}} Suppose that $B$ is a Banach  $A-bimodule$ and $A$ has a $BAI$. Then we have the following assertions.
 \begin{enumerate}
\item  If $B^*$ factors on the left [resp. right] with respect to $A$ and $A$ has $Rw^*w-$ property [resp. $Lw^*w-$ property], then $Z_{B^{**}}(A^{**})=A^{**}$.
\item  If $B^*$ factors on the left [resp. right] with respect to $A$ and $A$ has $Rw^*w-$ property [resp. $Lw^*w-$ property] with respect $B$, then $Z_{A^{**}}(B^{**})=B^{**}$.
\end{enumerate}
\begin{proof}
 \begin{enumerate}
\item   Assume that $B^*$ factors on the left and $A$ has $Rw^*w-$ property. Let   $(b_{\alpha}^{\prime\prime})_{\alpha}\subseteq B^{**}$ such that  $b^{\prime\prime}_{\alpha} \stackrel{w^*} {\rightarrow}b^{\prime\prime}$. Since $B^*A=B^*$, for all $b^\prime\in B^*$ there are $x\in A$ and $y^\prime\in B^*$ such that $b^\prime=y^\prime x$. Then for all $a\in A$, we have
$$<\pi_\ell^{**}(b_\alpha^{\prime\prime},y^{\prime})x,a>=<b_\alpha^{\prime\prime},\pi_\ell^{*}(y^{\prime},a)x>
=<\pi_\ell^{**}(b_\alpha^{\prime\prime},b^{\prime}),a>$$$$=<b_\alpha^{\prime\prime},\pi_\ell^{*}(b^{\prime},a)>
\rightarrow<b^{\prime\prime},\pi_\ell^{*}(b^{\prime},a)>
=<\pi_\ell^{**}(b^{\prime\prime},y^{\prime})x,a>.$$
Thus, we conclude that $\pi_\ell^{**}(b_\alpha^{\prime\prime},y^{\prime})x\stackrel{w^*} {\rightarrow}
<\pi_\ell^{**}(b^{\prime\prime},y^{\prime})x$.
Since $A$ has $Rw^*w-$ property, $\pi_\ell^{**}(b_\alpha^{\prime\prime},y^{\prime})x\stackrel{w} {\rightarrow}
<\pi_\ell^{**}(b^{\prime\prime},y^{\prime})x$.
Now let $b^{\prime\prime}\in A^{**}$. Then
$$<\pi_\ell^{***}(a^{\prime\prime},b^{\prime\prime}_{\alpha}),b^\prime>=
<a^{\prime\prime},\pi_\ell^{**}(b^{\prime\prime}_{\alpha},b^\prime)>=
<a^{\prime\prime},\pi_\ell^{**}(b^{\prime\prime}_{\alpha},y^\prime)x>$$$$\rightarrow
<a^{\prime\prime},\pi_\ell^{**}(b^{\prime\prime},y^\prime)x>=
<\pi_\ell^{***}(a^{\prime\prime},b^{\prime\prime}),b^\prime>.$$
It follow that $a^{\prime\prime}\in Z_{B^{**}}(A^{**})=A^{**}$.\\
If $B^*$ factors on the right and $A$ has $Lw^*w-$ property, then proof is the same as preceding proof.\\
\item  Let $B^*$ factors on the left  with respect to $A$ and $A$ has $Rw^*w-$ property  with respect to $B$.
Assume that
$(a_{\alpha}^{\prime\prime})_{\alpha}\subseteq A^{**}$ such that  $a^{\prime\prime}_{\alpha} \stackrel{w^*} {\rightarrow}a^{\prime\prime}$. Since $B^*A=B$, for all $b^\prime\in B^*$ there are $x\in A$ and $y^\prime\in B^*$ such that $b^\prime=y^\prime x$. Then for all $b\in B$, we have
$$<\pi_r^{**}(a_\alpha^{\prime\prime},y^{\prime})x, b>=<\pi_r^{**}(a_\alpha^{\prime\prime},b^{\prime}), b>=
<a_\alpha^{\prime\prime},\pi_r^{*}(b^{\prime}, b)>$$$$=<a^{\prime\prime},\pi_r^{*}(b^{\prime}, b)>=<\pi_r^{**}(a^{\prime\prime},y^{\prime})x, b>.$$
Consequently  $\pi_r^{**}(a_\alpha^{\prime\prime},y^{\prime})x\stackrel{w^*} {\rightarrow}\pi_r^{**}(a^{\prime\prime},y^{\prime})x$. Since $A$ has $Rw^*w-$ property  with respect to $B$,
$\pi_r^{**}(a_\alpha^{\prime\prime},y^{\prime})x\stackrel{w} {\rightarrow}\pi_r^{**}(a^{\prime\prime},y^{\prime})x$.
It follow that for all $b^{\prime\prime}\in B^{**}$, we have
$$<\pi_r^{***}(b^{\prime\prime},a_\alpha^{\prime\prime}), b^\prime>=
<b^{\prime\prime},\pi_r^{**}(a_\alpha^{\prime\prime}, y^\prime)x>\rightarrow
<b^{\prime\prime},\pi_r^{**}(a^{\prime\prime}, y^\prime)x>$$
$$=<\pi_r^{***}(b^{\prime\prime},a^{\prime\prime}), b^\prime>.$$
Thus we conclude that $b^{\prime\prime}\in Z_{A^{**}}(B^{**})$.\\
The proof of the next assertions is the same as proceeding proof.
\end{enumerate}
\end{proof}

\noindent{\it{\bf Theorem 3-5.}} Suppose that $B$ is a Banach  $A-bimodule$. Then we have the following assertions.
\begin{enumerate}
\item  If $a_0\in A$ has $Rw^*w-$ property with respect to $B$, then $a_0A^{**}\subseteq Z_{B^{**}}(A^{**})$ and $a_0B^*\subseteq wap_\ell(B)$.
\item  If $a_0\in A$ has $Lw^*w-$ property with respect to $B$, then $A^{**}a_0\subseteq Z_{B^{**}}(A^{**})$ and $B^*a_0\subseteq wap_\ell(B)$.
\item  If $a_0\in A$ has $Rw^*w-$ property with respect to $B$, then $a_0B^{**}\subseteq Z_{A^{**}}(B^{**})$ and $B^*a_0\subseteq wap_r(B)$.
\item  If $a_0\in A$ has $Lw^*w-$ property with respect to $B$, then $B^{**}a_0\subseteq Z_{A^{**}}(B^{**})$ and $a_0B^*\subseteq wap_r(B)$.
\end{enumerate}
\begin{proof}
\begin{enumerate}

\item Let   $(b_{\alpha}^{\prime\prime})_{\alpha}\subseteq B^{**}$ such that  $b^{\prime\prime}_{\alpha} \stackrel{w^*} {\rightarrow}b^{\prime\prime}$. Then for all $a\in A$ and  $b^\prime\in B^*$, we have
$$<\pi_\ell^{**}(b_\alpha^{\prime\prime},b^{\prime})a_0,a>=<\pi_\ell^{**}(b_\alpha^{\prime\prime},b^{\prime}),a_0a>
=<b_\alpha^{\prime\prime},\pi_\ell^{*}(b^{\prime},a_0a)>$$
$$\rightarrow <b^{\prime\prime},\pi_\ell^{*}(b^{\prime},a_0a)>=<\pi_\ell^{**}(b^{\prime\prime},b^{\prime})a_0,a>.$$
It follow that $\pi_\ell^{**}(b_\alpha^{\prime\prime},b^{\prime})a_0\stackrel{w^*} {\rightarrow}
\pi_\ell^{**}(b^{\prime\prime},b^{\prime})a_0$. Since $a_0$ has $Rw^*w-$ property with respect to $B$,
$\pi_\ell^{**}(b_\alpha^{\prime\prime},b^{\prime})a_0\stackrel{w} {\rightarrow}
\pi_\ell^{**}(b^{\prime\prime},b^{\prime})a_0$.\\
We conclude that $a_0a^{\prime\prime}\in  Z_{B^{**}}(A^{**})$ so that   $a_0 A^{**}\in  Z_{B^{**}}(A^{**})$.
Since $\pi_\ell^{**}(b^{\prime\prime},b^{\prime})a_0=\pi_\ell^{**}(b^{\prime\prime},b^{\prime}a_0)$, $a_0B^*\subseteq wap_\ell(B)$.\\
\item proof is similar to (1).
\item Assume that
$(a_{\alpha}^{\prime\prime})_{\alpha}\subseteq A^{**}$ such that  $a^{\prime\prime}_{\alpha} \stackrel{w^*} {\rightarrow}a^{\prime\prime}$. Let $b\in B$ and $b^\prime\in B^*$. Then we have
$$<\pi_r^{**}(a_\alpha^{\prime\prime},b^{\prime})a_0, b>=<\pi_r^{**}(a_\alpha^{\prime\prime},b^{\prime}),a_0 b>=
<a_\alpha^{\prime\prime},\pi_r^{*}(b^{\prime},a_0 b)>$$
$$\rightarrow <a^{\prime\prime},\pi_r^{*}(b^{\prime},a_0 b)>=<\pi_r^{**}(a^{\prime\prime},b^{\prime})a_0, b>.$$
Thus we conclude $\pi_r^{**}(a_\alpha^{\prime\prime},b^{\prime})a_0\stackrel{w^*} {\rightarrow}\pi_r^{**}(a^{\prime\prime},b^{\prime})a_0$.
Since $a_0$ has $Rw^*w-$ property with respect to $B$, $\pi_r^{**}(a_\alpha^{\prime\prime},b^{\prime})a_0\stackrel{w} {\rightarrow}\pi_r^{**}(a^{\prime\prime},b^{\prime})a_0$.
If $b^{\prime\prime}\in B^{**}$, then we have
$$<\pi_r^{***}(a_0b^{\prime\prime},a_\alpha^{\prime\prime}), b^\prime>=<a_0b^{\prime\prime},\pi_r^{***}(a_\alpha^{\prime\prime}, b^\prime)>=
<b^{\prime\prime},\pi_r^{**}(a_\alpha^{\prime\prime}, b^\prime)a_0>$$
$$=<b^{\prime\prime},\pi_r^{**}(a_\alpha^{\prime\prime}, b^\prime)a_0>=<\pi_r^{***}(a_0b^{\prime\prime},a^{\prime\prime}), b^\prime>.$$
It follow that $a_0b^{\prime\prime}\in Z_{A^{**}}(B^{**})$. Consequently we have  $a_0B^{**}\in Z_{A^{**}}(B^{**})$.\\
The proof of the next assertion is clear.\\
\item Proof is similar to (3).\\\\
\end{enumerate}
\end{proof}

\noindent{\it{\bf Theorem 3-6.}} Let $B$ be a Banach  $A-bimodule$. Then we have the following assertions.
\begin{enumerate}
\item  Suppose
$$\lim_{\alpha}\lim_{\beta}<b^\prime_{\beta},b_{\alpha}>=\lim_{\beta}\lim_{\alpha}
<b^\prime_{\beta},b_{\alpha}>,$$ for every $(b_{\alpha})_{\alpha}\subseteq B$ and $(b^\prime_{\beta})_{\beta}\subseteq B^*$. Then $A$ has $Lw^*w-$ property and $Rw^*w-$ property with respect to $B$.\\
\item  If for some $a\in A$, $$\lim_{\alpha}\lim_{\beta}<ab^\prime_{\beta},b_{\alpha}>=\lim_{\beta}\lim_{\alpha}
<ab^\prime_{\beta},b_{\alpha}>,$$
 for every $(b_{\alpha})_{\alpha}\subseteq B$ and $(b^\prime_{\beta})_{\beta}\subseteq B^*$, then $a$ has $Rw^*w-$ property with respect to $B$. Also if for some $a\in A$, $$\lim_{\alpha}\lim_{\beta}<b^\prime_{\beta}a,b_{\alpha}>=\lim_{\beta}\lim_{\alpha}
<b^\prime_{\beta}a,b_{\alpha}>,$$
for every $(b_{\alpha})_{\alpha}\subseteq B$ and $(b^\prime_{\beta})_{\beta}\subseteq B^*$, then $a$ has $Lw^*w-$ property with respect to $B$.
\end{enumerate}
\begin{proof}
\begin{enumerate}
 \item  Assume that $a\in A$ such that $ab_{\beta}^\prime\stackrel{w^*} {\rightarrow}0$ where
$(b^\prime_{\beta})_{\beta}\subseteq B^*$. Let $b^{\prime\prime}\in B^{**}$ and $(b_{\alpha})_{\alpha}\subseteq B$
such that $b_{\alpha}\stackrel{w^*} {\rightarrow}b^{\prime\prime}$. Then
$$\lim_{\beta}<b^{\prime\prime},ab_{\beta}^\prime>=\lim_{\beta}\lim_{\alpha}<b_{\alpha},ab_{\beta}^\prime>=
\lim_{\beta}\lim_{\alpha}<ab_{\beta}^\prime ,b_{\alpha}>$$
$$=\lim_{\alpha}\lim_{\beta}<ab_{\beta}^\prime ,b_{\alpha}>=0.$$
We conclude that $ab_{\beta}^\prime\stackrel{w} {\rightarrow}0$, so $A$ has $Lw^*w-$ property. It also easy that $A$ has  $Rw^*w-$ property.\\
\item  Proof is easy and is the same as (1).
\end{enumerate}
\end{proof}
\noindent{\it{\bf Definition 3-7.}} Let $B$ be a left Banach $A-module$. We say that $B^*$ strong factors on the left [resp. right] if for all $(b^\prime_\alpha)_\alpha\subseteq B^*$ there are $(a_\alpha)_\alpha\subseteq A$ and $b^\prime\in B^*$ such that $b^\prime_\alpha=b^\prime a_\alpha$ [resp. $b^\prime_\alpha= a_\alpha b^\prime$] where $(a_\alpha)_\alpha$ has limit the $weak^*$ topology in $A^{**}$.\\
If $B^*$ strong factors on the left and right, then we say that $B^*$ strong factors on the both side.\\
It is clear that if $B^*$ strong factors on the left [resp. right], then $B^*$ factors on the left [resp. right].\\\\
\noindent{\it{\bf Theorem 3-8.}} Suppose that $B$ is a Banach  $A-bimodule$. Assume that $AB^*\subseteq wap_\ell B$. If
$B^*$ strong factors on the left [resp. right], then $A$ has $Lw^*w-$ property [resp. $Rw^*w-$ property ] with respect to $B$.
\begin{proof}  Let   $(b_{\alpha}^{\prime})_{\alpha}\subseteq B^{*}$ such that  $ab_\alpha^\prime\stackrel{w^*} {\rightarrow}0$. Since $B^*$ strong factors on the left, there are
$(a_\alpha)_\alpha\subseteq A$ and $b^\prime\in B^*$ such that $b^\prime_\alpha=b^\prime a_\alpha$. Let
$b^{\prime\prime}\in B^{**}$ and $(b_{\beta})_{\beta}\subseteq B$ such that
 $b_{\beta}\stackrel{w^*} {\rightarrow}b^{\prime\prime}$. Then we have
$$\lim_\alpha<b^{\prime\prime},ab_\alpha^\prime>=\lim_\alpha\lim_{\beta}<b_{\beta},ab_\alpha^\prime>
=\lim_\alpha\lim_{\beta}<ab_\alpha^\prime,b_{\beta}>$$
$$=\lim_\alpha\lim_{\beta}<ab^\prime a_\alpha,b_{\beta}>=\lim_\alpha\lim_{\beta}<ab^\prime, a_\alpha b_{\beta}>$$
$$=\lim_{\beta}\lim_\alpha<ab^\prime, a_\alpha b_{\beta}>=\lim_{\beta}\lim_\alpha<ab_\alpha^\prime,b_{\beta}>=0$$
It follow that $ab_\alpha^\prime\stackrel{w} {\rightarrow}0$.
\end{proof}
\noindent{\it{\bf Problems .}}
\begin{enumerate}
 \item Suppose that $B$ is a Banach  $A-bimodule$. If $B$ is left or right factors with respect to $A$, dose $A$ has $Lw^*w-$property or $Rw^*w-$property, respectively?
\item Suppose that $B$ is a Banach  $A-bimodule$. Let $A$ has $Lw^*w-$property with respect to $B$. Dose $Z_{B^{**}}(A^{**})=A^{**}$?
\end{enumerate}

\bibliographystyle{amsplain}

\begin{thebibliography}{99}
\bibitem{1} R. E.  Arens, {\it The adjoint of a bilinear operation}, Proc. Amer. Math. Soc. {\bf 2} (1951), 839-848.
\bibitem{2} N. Arikan, {\it A simple condition ensuring the Arens
regularity of bilinear mappings}, Proc. Amer. Math. Soc. {\bf 84}
(4) (1982), 525-532.
\bibitem{3} J. Baker, A.T. Lau, J.S. Pym {\it Module homomorphism and topological centers associated with
 weakly sequentially compact Banach algebras}, Journal of Functional Analysis. {\bf 158} (1998), 186-208.
\bibitem{4} F. F. Bonsall, J. Duncan, {\it Complete normed algebras}, Springer-Verlag, Berlin 1973.
\bibitem{5} H. G. Dales, A. Rodrigues-Palacios, M.V. Velasco, {\it The second transpose of a derivation}, J. London. Math. Soc. {\bf2} 64 (2001) 707-721.
  \bibitem{6} H. G. Dales, {\it Banach algebra and automatic continuity}, Oxford 2000.

\bibitem{7} N. Dunford, J. T. Schwartz, {\it Linear operators.I},
Wiley, New york 1958.
\bibitem{8} M. Eshaghi Gordji, M. Filali, {\it Arens regularity of module actions},  Studia Math. {\bf 181} 3 (2007), 237-254.
\bibitem{8} M. Eshaghi Gordji, M. Filali, {\it Weak amenability of the second dual of a Banach algebra}, Studia Math. {\bf182} 3 (2007), 205-213.




\bibitem{8} K. Haghnejad Azar, A. Riazi,  {\it Arens regularity of bilinear forms and unital Banach module space}, (Submitted).


\bibitem{9} E. Hewitt, K. A. Ross,  {\it Abstract harmonic analysis}, Springer, Berlin, Vol I 1963.
\bibitem{10} E. Hewitt, K.  A. Ross,  {\it Abstract harmonic analysis}, Springer, Berlin, Vol II 1970.
\bibitem{11}  A. T. Lau, V. Losert, {\it On the second Conjugate Algebra of locally
compact groups}, J. London Math. Soc.  {\bf 37} (2)(1988),
464-480.
\bibitem{12} A. T. Lau, A. $\ddot{U}lger$, {\it Topological center of certain dual
algebras}, Trans. Amer.  Math. Soc. {\bf 348} (1996), 1191-1212.

\bibitem{13} S. Mohamadzadih, H. R. E. Vishki, {\it Arens regularity of module actions and the second adjoint of a derivation}, Bulletin of the Australian Mathematical Society {\bf77} (2008), 465-476.
\bibitem{15} M. Neufang, {\it Solution to a conjecture by Hofmeier-Wittstock},
Journal of Functional Analysis. {\bf 217} (2004), 171-180.
\bibitem{16} M. Neufang, {\it On a conjecture
 by Ghahramani-Lau and related problem concerning topological center}, Journal of Functional Analysis.
  {\bf 224} (2005), 217-229.
  \bibitem{17}  J. S. Pym, {\it The convolution of functionals on spaces of bounded functions},
 Proc. London Math Soc.  {\bf 15} (1965), 84-104.
\bibitem{18}A. $\ddot{U}lger$, {\it Arens regularity sometimes implies the RNP}, Pacific Journal of
Math. {\bf 143} (1990), 377-399.
\bibitem{19} A. $\ddot{U}lger$, {\it Some stability properties of Arens regular bilinear operators}, Proc. Amer. Math. Soc. (1991) {\bf 34}, 443-454.
\bibitem{20}  A. $\ddot{U}lger$, {\it Arens regularity of weakly sequentialy compact Banach algebras},
Proc. Amer. Math. Soc. {\bf 127} (11) (1999), 3221-3227.
\bibitem{17} A. $\ddot{U}lger$, {\it Arens regularity of the algebra $A\hat{\otimes}B$}, Trans. Amer. Math. Soc. {\bf 305} (2) (1988) 623-639.
\bibitem{21} P. K. Wong, {\it The second conjugate algebras of
Banach algebras}, J. Math. Sci. {\bf 17} (1) (1994), 15-18.
 \bibitem{22} N. J. Young  {\it Theirregularity of multiplication in group algebra} Quart. J. Math. Soc., {\bf 24} (2)  (1973), 59-62.
\bibitem{22} Y. Zhing,  {\it Weak amenability of module extentions of Banach algebras} Trans. Amer. Math. Soc. {\bf 354} (10) (2002), 4131-4151.

\end{thebibliography}

\it{Department of Mathematics, Amirkabir University of Technology, Tehran, Iran\\
{\it Email address:} haghnejad@aut.ac.ir

\end{document}

%------------------------------------------------------------------------------
% End of journal.tex
%------------------------------------------------------------------------------